\documentclass[a4paper,12pt]{article}

\begin{document}

\thispagestyle{empty}

\vspace*{2cm}
\begin{center}
{\Large \bf Twist of Lie algebras by 6 dimensional subalgebra}\footnote{
Talk given at 32nd Symposium on Mathematical Physics, 
June 6-10, 2000, Torun, Poland}

\vspace{1cm}
N. Aizawa\footnote{
On leave of absence 
from Department of Applied Mathematics, 
Osaka Women's University, Sakai, Osaka 590-0035, Japan}

Arnold Sommerfeld Institute for Mathematical Physics, \\
Technical University of Clausthal,  \\
Leibnizstra\ss e 10, D-38678 Claustahl-Zellerfeld, Germany

\end{center}

\vfill
\begin{abstract}
  A new non-standard deformation of all types of classical Lie algebras 
is constructed by means of Drinfel'd twist based on a six dimensional 
subalgebra. This is an extension of extended twists introduced by 
Kulish $et\ al.$ For the algebra $ {\cal M}_3 \simeq so(3,2) $, 
a relation to a known non-standard deformation is discussed.
\end{abstract}

\newpage
\section{Introduction}

  It is known that quantum algebras belong to one of two types of 
Hopf algebras: 
quasitriangular and triangular \cite{qg}. The $q$-deformed algebras 
by Drinfel'd and Jimbo belong to the first type. The algebras belong 
to the second type are called Jordanian quantum algebras or 
non-standard quantum algebras. The typical example of Jordanian 
quantum algebras is a deformation of $ sl(2) $ introduced by 
Ohn \cite{ohn}. In general, triangular quantum algebras are obtained 
by applying Drifel'd twist to Lie algebras \cite{twist}. The twisting 
that produces a isomorphic algebra to Ohn's one is found in 
\cite{jor,jor2}.

  In this article, we discuss a new non-standard deformation of all types 
of classical Lie algebras, that can be regarded as an extension of 
non-standard quantum algebras introduced recently by Kulish $et\ al.$ 
\cite{et}. To make clear the 
novelty of our results and to fix notations, we briefly recall the 
definition of twisting and already known twisting for Lie algebras in the 
next section. The motivation of this study will be mentioned there. 
\S 3 is the main part of this article where the new non-standard 
deformation is introduced. As an physical example, we take an algebra 
of conformal transformations in (2+1) dimensional Minkowskian spacetime 
in \S 4. The equivalence of our result and another non-standard 
deformation by Herranz \cite{her} is discussed. \S 5 is devoted to 
concluding remarks.

  Drinfel'd developed the idea of twisting in his study of quasi-Hopf algebras, 
however, we restrict ourselves to ordinary Hopf algebras throughout this 
article.

\section{Brief review of Drinfel'd twist}

  Let $ ({\cal H}_0, m_0, \Delta_0, \epsilon_0, S_0) $ be a Hopf algebra, 
$ m_0, \Delta_0, \epsilon_0 $ and $ S_0 $ denote product, coproduct, 
counit and antipode of $ {\cal H}_0$, respectively. 
Suppose that an invertible element 
$ {\cal F} \in {\cal H}_0 \otimes {\cal H}_0 $, called twistor or 
twist element, satisfies
\begin{equation}
  {\cal F}_{12} (\Delta_0 \otimes id)({\cal F}) = 
  {\cal F}_{23} (id \otimes \Delta_0)({\cal F}),
  \label{twist1}
\end{equation}
\begin{equation}
  (\epsilon_0 \otimes id)({\cal F}) = (id \otimes \epsilon_0)({\cal F}) = 1,
  \label{twist2}
\end{equation}
where $ {\cal F}_{12} = {\cal F} \otimes 1,\ 
{\cal F}_{23} = 1 \otimes {\cal F} $ and $id $ denotes the identity 
mapping as usual. Then $ ({\cal H}_0, m_0, \Delta, \epsilon_0, S) $ with
\begin{equation}
 \Delta = {\cal F} \Delta_0 {\cal F}^{-1},  \quad
  S = v S_0 v^{-1}, \label{twistedco}
\end{equation}
is also a Hopf algebra. We denote this new Hopf algebra by $ {\cal H} $. 
The element $ v \in {\cal H}_0 $ is obtained from the twistor. 
Writing $ {\cal F} = \sum_i f_i \otimes f^i $, the element $ v$ is given by 
$ v = \sum_i f_i S_0(f^i). $ Note that, in the new Hopf algebra ${\cal H}$, 
only coproduct and antipode are twisted, while product and coproduct are not 
changed. The universal $R$-matrix for $ {\cal H} $ is also obtained from the 
one for ${\cal H}_0$,
\begin{equation}
  {\cal R} = {\cal F}_{21} R_0 {\cal F}^{-1}, \label{RR0}
\end{equation}
where $ {\cal R} $ and $ {\cal R}_0 $ are the universal $R$-matrices for 
$ {\cal H} $ and ${\cal H}_0$, respectively. 

  In the case of non-standard quantum algebras, ${\cal H}_0$ is a Lie 
algebra (strictly speaking, universal enveloping algebra of a Lie algebra). 
Since the twisting does not change the product, non-standard quantum 
algebras have the same commutation relations as Lie algebras. However, 
they have deformed coproducts and antipodes. The triangularity 
of non-standard algebras stems from the fact that ${\cal R}_0 $ for 
Lie algebras is given by $ 1 \otimes 1$. 
The advantages of non-standard quantum 
algebras are: (i) since they have undeformed commutation relations, 
we know their irreducible representations, (ii) explicit form of 
universal $R$-matrix is known. Combining (i) and (ii), we obtain matrix 
representations of universal $R$-matrix, this gives us further advantages: 
(iii) dual quantum groups are easily obtained by 
FRT-formalism \cite{frt}, (iv) covariant differential calculus are easily 
obtained, and so on.

  As is clear from the above discussion, once we obtain a explicit form 
of twistor, the rest steps of construction of non-standard quantum 
algebras and groups are automatically carried out. Note that the possibility 
of twisting for a Lie algebra is not unique, namely, a Lie algebra admits 
some twistors. Note also that non-standard algebras have some physical 
applications (see $e.g.$ \cite{ks}, \cite{bhp} and \cite{bch}) 
and quantum groups are convenient tool 
to describe non-commutative spacetime. Therefore, we think that it is 
important to investigate possible twistors for Lie algebras. 
In the study of twisting, the following observation is useful: 
Let $ {\cal A} $ be 
a subalgebra of a Lie algebra ${\cal H}_0$ and $ {\cal F}$ is a twistor 
for the subalgebra $ {\cal A}$. Then, the twistor ${\cal F}$ can be 
regarded as a twistor for $ {\cal H}_0$ and we can twist $ {\cal H}_0 $ 
with this twistor. In the following, 
we give five examples of twistors classified according 
to the subalgebra, then we discuss an extension of one of them. 

  Probably, the most wellknown twistors are the socalled Reshetikhin 
twist \cite{res} and Jordanian twist \cite{jor,jor2}. 
In the Reshetikhin twist, the Cartan subalgebra is used as 
the subalgebra. Thus any Lie algebras of rank $ \geq 2$ 
admit Reshetikhin twist. 
On the other hand, the subalgebra for Jordanian twist 
is the Borel subalgebra: $ \{ H,\; E \ | \ [H,\; E] = 2E \}. $ 
The explicit form of twistor is given by
\begin{equation}
  {\cal F}_J = \exp( -\frac{1}{2} H \otimes \sigma ),  \quad 
  \sigma = - \ln(1 - zE), \label{Fj}
\end{equation}
where $z$ is a deformation parameter. We denote a deformation parameter 
by $z$ throughout this article and $ z = 0$ corresponds to undeformed 
limit. Three nontrivial extensions of Jordanian twist were considered 
very recently \cite{et,pet,chain}. 
In \cite{et}, the Borel subalgebra is extended to four dimensional one 
that is a semidirect sum of Borel subalgebra and two additional elements 
$ A, B. $ Let $ {\cal A}_E = \{H, E, A, B \} $ be the subalgebra subject 
to the relations
\begin{eqnarray}
     & & [H, E] = \delta E, \quad [H, A] = \alpha A, \quad [H, B] = \beta B, 
         \label{subet} \\
     & & [A, B] = \gamma E, \quad [E, A] = [E, B] = 0, \quad
     \alpha + \beta = \delta. \nonumber
\end{eqnarray}
Then the following is a twistor
\begin{equation}
   {\cal F}_E = \exp(A \otimes B e^{-\beta \sigma / \delta} ) {\cal F}_J.
   \label{etwistor}
\end{equation}
This twisting is called extended (Jordanian) twist. 
It is verified that all types of classical Lie algebra have the 
subalgebra $ {\cal A}_E. $  We use different convention from \cite{et} 
for the definition of $\sigma$, but 
this may not cause any confusion. The extension considered in 
\cite{pet} is certain limits of extended twists and called peripheric 
extended twists. The extended twists have nontrivial limits for 
$ \alpha \rightarrow 0 $ or $ \beta \rightarrow 0.$ 
One can verify that the peripheric extended twists are applicable to 
inhomogeneous Lie algebras $ isu(n) $ and $ iso(n). $ 
In \cite{chain}, regular injections 
$ {\cal A}_p \subset {\cal A}_{p-1} \subset \cdots \subset 
{\cal A}_1 \subset {\cal A}_0 $ of Lie algebras and twistors of 
extended twists are considered. It is shown that a product of twistors of 
extended twists corresponding to each subsets 
$ {\cal A}_k (k = 0, 1, \cdots, p) $ can produce a new twistor for 
the following sequences:
\begin{eqnarray}
  & & sl(n) \supset sl(n-2) \supset \cdots \supset sl(n-2k) \supset \cdots,
      \nonumber \\
  & & so(2n) \supset so(2n-4) \supset \cdots \supset so(2n-4k) \supset \cdots,
      \nonumber \\
  & & so(2n+1) \supset \cdots \supset so(2n-4k+1) \supset \cdots
      \nonumber
\end{eqnarray}
The case of $ sp(n) $ is also considered in \cite{chain}, though the situation 
is different from other classical Lie algebras. 

  All twistors mentioned above except Rshetikhin twist have common 
properties: (i) twisted coproduct for $ \sigma $ 
is primitive, that is, 
$ \Delta(\sigma) = \sigma \otimes 1 + 1 \otimes \sigma, $ (ii) the 
twistors are factorizable, that is, they satisfy the relations
\begin{equation}
(\Delta_0 \otimes id) ({\cal F}) = {\cal F}_{13} {\cal F}_{23}, \quad
 (id \otimes \Delta) ({\cal F}) = {\cal F}_{12} {\cal F}_{13}.
 \label{factorize}
\end{equation}
These relations guarantee that the $ {\cal F} $ satisfies the 
condition (\ref{twist1}).

  In the next section, we consider a further extension of extended twists in 
such a way that the four dimensional subalgebra $ {\cal A}_E $ is 
replaced with six dimensional one.

\section{New twisting for classical Lie algebras}

  Let us consider an algebra $ {\cal A} $ of six elements 
$ H_i, E_i, A , B \ (i = 1, 2) $ satisfying
\begin{eqnarray}
  & & [H_i,\; E_i] = 2E_i,  \nonumber \\
  & & [H_1, \; H_2] = [E_1,\; E_2] = [H_1,\; E_2] = [H_2, \; E_1] = 0,
      \label{6dsub} \\
  & & [H_1,\; A] = -A, \quad [H_1,\; B] = B, \quad [H_2,\; A] = A, 
      [H_2,\; B] = B,  \nonumber \\
  & & [A, \; E_1 ] = 2B, \quad [A,\; E_2 ] = 0, \quad
      [E_i,\; B] = 0, \quad [A, B] = E_2. \nonumber
\end{eqnarray}
The four elements $ \{H_2, E_2, A, B \} $ form the subalgebra 
$ {\cal A}_E $ of extended twists ($\alpha=\beta=\gamma=1 $). 
The additional elements $ H_1, E_1 $ 
form a Borel subalgebra, thus the algebra $ {\cal A} $ is a 
semidirect sum of $ {\cal A}_E $ and an extra Borel subalgebra.
We can also regard the algebra $ {\cal A} $ as a semidirect sum of 
$ \{A, B\} $ and a direct sum of two Borel subalgebras 
\[
 {\cal A} = (\{H_1, E_1\} \oplus \{H_2, E_2\})  
  \in\hspace{-.67em}{\scriptstyle \mid} \ \{A, B \}.
\] 
The following invertible element $ {\cal F} $ satisfies the definition 
of twistor.
\begin{equation}
 {\cal F} = \exp(-\frac{1}{2} H_1 \otimes \sigma_1) 
       \exp(-z A \otimes B e^{\frac{1}{2} \sigma_2}) 
       \exp(-\frac{1}{2} H_2 \otimes \sigma_2)_,
    \label{6dtwistor}
\end{equation}
where
\begin{equation}
  \sigma_1 = - \ln( 1 - z (E_1 + z B^2 e^{\sigma_2})), \quad 
  \sigma_2 = - \ln( 1 - z E_2 ).
\end{equation}
The two factors from the right are a extended twist and left most factor 
does not commute with the rest part of ${\cal F}$. Therefore this is a 
nontrivial extension of extended twists. The twistor (\ref{6dtwistor}) has 
similar properties as extended twists. Namely, the twisted coproducts for 
$ \sigma_i $ are primitive
\begin{equation}
 \Delta(\sigma_i) = \sigma_i \otimes 1 + 1 \otimes \sigma_i, \quad 
  i = 1, 2,
 \label{coprosigma}
\end{equation}
and the ${\cal F}$ is factorizable (see eq.(\ref{factorize})).

  To prove the above statements, we first calculate the twisted 
coproducts of the elements of ${\cal A}$ by the twistor (\ref{6dtwistor}). 
They are given by
\begin{eqnarray}
 & & \Delta(H_1) = H_1 \otimes e^{\sigma_1} + 1 \otimes H_1, \nonumber \\
 & & \Delta(E_1) = E_1 \otimes e^{-\sigma_1} + 1 \otimes E_1 
    - 2z B \otimes \rho e^{-(\sigma_1 + \sigma_2)/2} 
    + z^2 E_2 \otimes \rho^2 e^{-\sigma_2}_, \nonumber \\
 & & \Delta(H_2) = H_2 \otimes e^{\sigma_2} + 1 \otimes H_2 
    + 2z A \otimes \rho e^{(\sigma_1 + \sigma_2)/2} 
    + z^2 H_1 \otimes \rho^2 e^{\sigma_1}_, \label{copro} \\
 & & \Delta(E_2) = E_2 \otimes e^{-\sigma_2} + 1 \otimes E_2, \nonumber \\
 & & \Delta(A) = A \otimes e^{(\sigma_1 + \sigma_2)/2} + 1 \otimes A
     + z H_1 \otimes (B + z E_2 \rho) e^{\sigma_1}_, \nonumber \\
 & & \Delta(B) = B \otimes e^{-(\sigma_1 + \sigma_2)/2} 
      + e^{-\sigma_2} \otimes B, \nonumber
\end{eqnarray}
where $ \rho = B e^{\sigma_2}. $ With these coproducts, we can verify that 
the $ \sigma_i$'s are primitive. 
It follows that the twistor (\ref{6dtwistor}) satisfies the factorizable 
relations (\ref{factorize}). Thus we have proved that the ${\cal F}$ 
satisfies the condition (\ref{twist1}). The condition (\ref{twist2}) is 
easily verified by noticing that $ \epsilon_0(X) = 0 $ for all elements 
of any Lie algebras.

  The universal $R$-matrix of the twisted Lie algebra is given by
\begin{eqnarray}
  R ={\cal F}_{21} {\cal F}^{-1} 
    &=& \exp(-\frac{1}{2} \sigma_1 \otimes H_1) \exp(-z \rho \otimes A)
        \exp(-\frac{1}{2} \sigma_2 \otimes H_2)
    \nonumber \\
    &\times & \exp(\frac{1}{2}H_2 \otimes \sigma_2) 
        \exp(z A \otimes \rho)
        \exp(\frac{1}{2} H_1 \otimes \sigma_1). \label{quantumr}
\end{eqnarray}
One can easily write down the corresponding classical $r$-matrix 
that solves the classical Yang-Baxter equation by keeping up to the 
first order of the deformation parameter
\begin{equation}
 r = \frac{1}{2} H_1 \wedge E_1 + \frac{1}{2} H_2 \wedge E_2 + A \wedge B.
 \label{classcalr}
\end{equation}

 We next show that all types of classical Lie algebras have the six 
dimensional subalgebra $ {\cal A}, $ namely, we obtain new 
non-standard quantum algebras. For $ sl(n), $ it is convenient to 
work with canonical basis
\begin{equation}
 [E_{ab},\; E_{cd}] = E_{ad} \delta_{bc} - E_{cb} \delta_{ad}, \quad 
 a, b, c, d = 1, \cdots, n.
 \label{cansln}
\end{equation}
In this case, the six dimensional subalgebra $ {\cal A} $ is found for 
$ n \geq 4 $
\begin{eqnarray}
  & & H_1 = \sum_{n/2 \geq k \geq 2} (E_{kk} - E_{n-k+1, n-k+1}), \quad 
      E_1 = \sum_{n/2 \geq k \geq 2} E_{k, n-k+1}, \nonumber \\
  & & H_2 = E_{11} - E_{nn}, \quad E_2 = E_{1n}, \label{forsln} \\
  & & A = 2 \sum_{n/2 \geq k \geq 2} b^{1,n-k+1} E_{1k} 
        - 2 \sum_{n-1 \geq \lambda > n/2} b^{n-\lambda+1,n} E_{\lambda n},
      \nonumber \\
  & & B = \sum_{n/2 \geq k \geq 2} b^{k,n} E_{kn} 
        + \sum_{n-1 \geq \lambda > n/2} b^{1,\lambda} E_{1\lambda}, 
  \nonumber
\end{eqnarray}
where the complex coefficients $ b^{a,b}$'s have to satisfy
\begin{equation}
 4 \sum_{n/2 \geq k \geq 2} b^{1,n-k+1} b^{k,n} = 1.
 \label{condsln}
\end{equation}
This condition on the coefficients stems from the commutation relation 
$ [A,\; B] = E_2. $ Other commutation relations are hold for any values of 
$ b^{a,b}_. $

  We also use the canonical basis for $so(n)$. 
\begin{eqnarray}
   & & [Y_{ab},\; Y_{cd} ] = i (Y_{ad} \delta_{bc} + Y_{bc} \delta_{ad} 
      - Y_{ac} \delta_{bd} - Y_{bd} \delta_{ac}),
    \label{canson} \\
   & &  Y_{ab} = - Y_{ba}, \quad a, b, c, d = 1, \cdots, n. \nonumber
  \label{defson}
\end{eqnarray}
In this case, the $ {\cal A} $ is found for $ n \geq 5. $
\begin{eqnarray}
  & & H_1 = Y_{1n} + Y_{1\; n-1}, \quad
      E_1 = Y_{1\; n-1} - Y_{2n} + iY_{12} + iY_{n-1\; n}, \nonumber \\
  & & H_2 = Y_{1n} - Y_{1\; n-1},\quad
      E_2 = \frac{1}{2}(Y_{1\; n-1} + Y_{2n} -i Y_{12} + i Y_{n-1\; n}),
      \label{forson} \\
  & & A = \sum_{k=3}^{n-2} a^k (Y_{n-1\; k} -i Y_{2k}), \quad 
      B = \sum_{k=3}^{n-2} a^k (Y_{kn} - iY_{1k}).
      \nonumber
\end{eqnarray}
The commutation relation $ [A,\; B] = E_2 $ imposes a condition 
on the coefficients $a^k_.$
\begin{equation}
  2 \sum_{k=3}^{n-2} (a^k)^2 = 1, \qquad a^k \in {\rm \bf C}. 
  \label{conson}
\end{equation}

  For $ sp(2n), $ the six dimensional subalgebra $ {\cal A} $ is found for 
$ n \geq 2. $ In terms of the canonical basis
\begin{eqnarray}
 & & [Z_{ab},\; Z_{cd}] = {\rm sign}(bc) 
      ( Z_{ad} \delta_{bc} + Z_{-b\; -c} \delta_{ad}
     + Z_{a\; -c} \delta_{-b\; d} + Z_{-b\; d} \delta_{c\; -a} ),
     \label{defsp2n} \\
 & & Z_{ab} = -{\rm sign}(ab) Z_{-b\; -a}, \quad a, b = \pm1, \cdots, \pm n,
    \nonumber 
\end{eqnarray}
the elements of $ {\cal A} $ are given by
\begin{eqnarray}
  & & H_1 = \sum_{k=2}^{n} Z_{kk}, \quad 
      E_1 = \sum_{k=2}^n Z_{k\; -k}, \quad 
      H_2 = Z_{11}, \quad \ \; E_2 = Z_{1\; -1}, 
      \label{forsp2n} \\
  & & A = \sum_{k=2}^n a^k Z_{1k}, \quad 
      B = \sum_{k=2}^n a^k Z_{k\; -1}. \nonumber
\end{eqnarray}
A condition on the coefficients $ a^k$'s are obtained 
same way as $ sl(n) $ and $ so(n)$
\begin{equation}
 \sum_{k=2}^n (a^k)^2 = 1, \quad a^k \in {\rm \bf C}.
 \label{consp2n}
\end{equation}

  We have seen that all types of classical Lie algebras can be twisted 
by the six dimensional subalgebra ${\cal A}$. Other combinations of 
elements of Lie algebras could realize the subalgebra $ {\cal A}. $ 
An appropriate choice could be found when physical applications of 
twisted algebras are considered.

\section{Conformal algebra of (2+1)-dimensional spacetime}

  In this section, we apply the twisting by $ {\cal A} $ to the algebra 
$ {\cal M}_3 \simeq so(3,2) $: the algebra of conformal transformations 
in the (2+1)-dimensional Minkowskian spacetime. We consider the following 
basis of $ {\cal M}_3 = \{J, P_{\mu}, K_i, C_{\mu}, D \} $ where 
$ \mu = 0, 1, 2,\  i = 1, 2 $ and $J$ is a generator of rotations, 
$ P_0$ time translations, $ P_i $ space translations, $ K_i $ boosts, 
$ C_{\mu} $ special conformal transformations, and $ D$ dilatations. 
The commutation relations of the algebra are given in eq.(1.1) of 
\cite{her}. The subalgebra $ {\cal A} $ is given by
\begin{eqnarray}
  & & H_1 = D + K_1, \qquad E_1 = P_0 + P_1, \qquad 
  H_2 = D - K_1, \label{so32} \\
  & & E_2 = P_0 - P_1, \qquad A = K_2 + J, \qquad \ \ B = P_2. \nonumber
\end{eqnarray}
Therefore, we have obtained a twisted $ {\cal M}_3 $ that has undeformed 
commutation relations but deformed coproducts. 
It may be remarkable that no generators of conformal transformations 
appear in the algebra $ {\cal A}$. 

  On the other hand, another non-standard deformation of $ {\cal M}_3 $ is 
considered in \cite{her}. We denote it by $ U_z({\cal M}_3) $ and its 
generators by $ \tilde J, \tilde P_{\mu}\ etc$. 
The algebra $ U_z({\cal M}_3) $ is defined by deformed commutation 
relations and deformed coproducts (eqs.(3.3)-(3.6) of \cite{her}). 
The universal $R$-matrix for $ U_z({\cal M}_3) $ is not known. 
We here give two observations that imply an equivalence between our 
twisted $ {\cal M}_3 $ and $ U_z({\cal M}_3). $ First, there exist an 
invertible mapping between the Weyl subalgebras of twisted 
$ {\cal M}_3 $ and $ U_z({\cal M}_3) $
\begin{eqnarray}
 & & \tilde D = D, \quad \tilde J = J, \quad \tilde K_i = K_i,  
 \label{mapping} \\
 & & \tilde P_0 = \frac{1}{2z} (\sigma_1 + \sigma_2), \quad 
     \tilde P_1 = \frac{1}{2z} (\sigma_1 - \sigma_2), \quad 
     \tilde P_2 = P_2 e^{\sigma_2}. \nonumber
\end{eqnarray}
By this mapping, both undeformed commutation relations and deformed 
coproducts of twisted ${\cal M}_3 $ are transformed into the 
corresponding ones of $ U_z({\cal M}_3). $ The mapping for $ C_{\mu} $ 
is not known. Probably, it has quite complicated form. Second, 
twisted $ {\cal M}_3 $ and $ U_z({\cal M}_3) $ have the same classical 
$r$-matrix
\begin{equation}
  r = z (D \wedge P_0 + K_1 \wedge P_1 + (K_2+J) \wedge P_2).
  \label{classrm3}
\end{equation}
This classical $r$-matrix (and more general form with two 
deformation parameters) was firstly discussed in \cite{jmm}. 
If the equivalence is true,  by using (\ref{mapping}) we can easily 
write down the universal $R$-matrix for $ U_z({\cal M}_3) $ 
in terms of its generators. 

  Two contractions of $ {\cal M}_3 $, that produce conformal algebras of the 
(2+1) Galilean and the Carroll spacetime, are also considered in \cite{her} 
and corresponding contractions for $ U_z({\cal M}_3) $ are analyzed. It is 
impossible to apply the contractions to our twistor, since the contracted 
algebras do not have the six dimensional subalgebra $ {\cal A} $. 
However, we can find twistors for each contracted algebras and these twistors 
give the same classical $r$-matrices as \cite{her}. 
For the conformal algebra of the Galilean spacetime,
\begin{equation}
  {\cal F} = \exp(-zK_1 \otimes P_1) \exp(-zK_2 \otimes P_2).
  \label{twistorgal}
\end{equation}
The commutation relations of this algebra are given in eq.(1.4) of \cite{her}. 
Two factors of (\ref{twistorgal}) commute each other. For the 
conformal algebra of the Carroll spacetime
\begin{equation}
   {\cal F} = \exp(-D \otimes \sigma) \exp(-z K_1 \otimes P_1 e^{\sigma}) 
         \exp(-z K_2 \otimes P_2 e^{\sigma}),
   \label{twistorcar}
\end{equation}
where $ \sigma = - \ln(1-zP_0). $ The commutation relations of this 
algebra are eq.(1.5) of \cite{her}. 
Three factors of (\ref{twistorcar}) commute one another.

\section{Concluding remarks}

  In this article, we have shown a new twistor that is an extension of 
the extended twists 
is applicable to all types of classical 
Lie algebras. Consequently new non-standard deformations of 
$ sl(n), so(n) $ and 
$ sp(2n) $ was obtained. When the twistor is applied to the algebra 
$ {\cal M}_3 $, the obtained algebra seems to be equivalent to a non-standard 
deformation of $ {\cal M}_3 $ by Herranz. 

  It is natural to ask whether the peripheric extended twists have a similar 
extension. Since the twists discussed in \S3 do not contain free parameters 
($ \alpha, \beta $ and $ \gamma $ of the extended twists), we can not 
repeat the same discussion as \cite{pet}. However, it turns out that the 
peripheric extended twists have an extension \cite{naru} where a five 
dimensional subalgeba is used instead of the six dimensional one. 
This extension of peripheric extended twists is appropriate to 
deform inhomegeneous Lie algebras. 

  Finally, we mention a physical application of the twisted algebras. 
It is known that symmetries of free massless Klein-Gordon equation are 
given by conformal algebras \cite{bx}. We mean, by symmetry, transformations 
of solutions into other solutions. If a difference analogue of 
Klein-Goldon equation on the uniform lattice is considered 
in (2+1)-dimensional spacetime, 
the twisted $ {\cal M}_3 $ can be regarded as its symmetry algebra. 
This will be discussed elsewhere. The similar situation was considered 
for (1+1)-dimensional free Schr\"odinger equation \cite{bhp,bch}. 

 After finishing this article, a paper \cite{defcar} that 
obtains the same result for $ so(5) $ was submitted to 
the preprint archive. 

\section*{Acknowledgments}

 The author would like to express his sincere thanks to 
Prof. H.-D. Doebner for his warm hospitality at Technical 
University of Clausthal. He also thanks Dr. F. J. Herranz for 
sending me some of his papers. He is grateful to Prof. 
Lukierski for pointing out the paper \cite{jmm} to him. 
The work was supported by 
the Ministry of Education, Science, Sports and Culture, Japan.

\end{document}